\documentclass[11pt]{amsart}

\usepackage{amsmath}
\usepackage{amssymb}
\usepackage{amsthm}
\usepackage{graphicx}

\newtheorem{theorem}{Theorem}[section]
\newtheorem*{theorem*}{Theorem}

\newtheorem*{lemma*}{Lemma}

\newtheorem*{criterion*}{Criterion}

\theoremstyle{definition}

\newcommand{\h}{\hspace{0.3cm}}

\numberwithin{equation}{section}

\title[Crest Factor and Turbulence in PDEs]{On an Analytical Criterion for detecting  Intermittent 
Turbulent Behaviour of Solutions of Partial Differential Equations}

\author[M. Bartuccelli]{Michele V Bartuccelli}

\address[M. V. Bartuccelli]{School of Mathematics and Physiscs, University of Surrey, Guildford, GU2 7XH, UK}
\email{{\tt m.bartuccelli@surrey.ac.uk}}

\author[G. Gentile]{Guido Gentile}
\address[G. Gentile]{Dipartimento di Matematica e Fisica, Universit\`a Roma Tre, Roma, I-00146,
Italy}
\email{\tt guido.gentile@uniroma3.it}



\begin{document}

\date{}

\maketitle

\vspace{-.2cm}
\centerline{\footnotesize Contribution to celebrate the 75th birthday of Professor Peter Constantin}
\vspace{-.2cm}

\begin{abstract}
 A main question in the study of partial differential equations 
is the following: how do we understand the nature of the solutions and, in particular,
how do we determine if a given solution  shows turbulent or non-turbulent behaviour?
Being able to answer such a question  would be a major advance in the comprehension of the nature of turbulence.
In this paper we focus on the case of intermittent turbulence and provide an analytical criterion,
based on the crest factor, which captures the essential feature of the solutions.
By computing the crest factor for the solutions of some classical equations, both linear and nonlinear,
we illustrate the capability of the criterion for discerning between solutions exhibiting time-intermittent turbulence behaviour
and solutions which either are not turbulent or show statistically stationary turbulence, like, for example, in the case described by Kolmogorov's theory.
\end{abstract}




\vspace{.3cm}

\section{\large Introduction} \label{sec1}

The American Nobel Prize laureate for Physics Richard Feynman
is reported to have designated turbulence as ``the most important unsolved problem of classical physics'',
because a description of the phenomenon  from first principles does not exist (a statement
along these lines is found in \cite[Vol.~I, \S3.7]{Feynman}).
In fact, turbulence is still regarded as one of the most important problems in 
mathematical physics today. The main difficulties encountered in trying to understand the nature 
of turbulence are related to its random behaviour both in space and time, due to the chaotic 
evolution of the fluid particles for large Reynolds numbers. Also there are many different length 
scales involved in the motion of the fluid as a whole which interact with one another in a very 
complicated way; this interaction makes it extremely difficult to predict the evolution of any 
(even smooth) initial condition for large Reynolds numbers and for large time.

There is, of course, a huge literature on the problem of turbulence and, therefore, we are forced 
to restrict ourselves to give some classical references without claiming to be complete or exhaustive:
for an introduction to the subject see, for example, 
\cite{CFT1,CFMT,CFT2,CO94,Fris,Foias1,Foias2,CGKMU,David,Berselli,Gal,CO2008,BaTi13}. 
A distinctive feature of turbulence, in all its manifestations, is chaoticity.
Kolmogorov's theory assumes the rate of energy transfer to change smoothly with space and time,
in both the quasi-linear regime (weak turbulence) and the fully nonlinear regime (strong turbulence).
On the other hand, experimental and numerical results show strong intermittent bursts
which give rise to non-uniform energy cascades and, as a consequence, to spiky, rare and irregular fluctuations (intermittency).
It is mainly the latter kind of hard turbulence that we have in mind in this paper.

Given the serious difficulties faced in investigating turbulence, it is apparent that it would be very
helpful to devise some analytical tools which shed some light on this subtle conundrum.
One of these tools is the so called crest factor \cite{MVB,MBPRS,BDG,BaGe}. 

Essentially, considering a solution of a PDE, the crest factor is the ratio between the $L^\infty$ and the $L^2$ norms
of the solution -- a more formal and detailed discussion is provided in Section \ref{SectionCF}.
In this paper, we argue that, for hard and time-intermittent turbulence to be exhibited by the solution, the corresponding crest factor must be large and,
besides depending explicitly on time, must be chaotic and characterised by strong fluctuations.
This provides a criterion to detect time-intermittent turbulent behaviour in a PDE. Of course, in the case of hard turbulence,
one expect to have strong intermittency in space as well. However, the criterion we propose allows us to recognise
immediately that intermittent turbulence is to be excluded if the crest factor is either time-independent or has a non-chaotic
dependence on time (for instance, if it is a periodic or quasi-periodic function). In addition, also in the case of non-intermittent turbulence,
from the size of the crest factor one can infer  what kind of turbulence may occur:
mild turbulence if the crest factor is of order one, hard turbulence if it is much larger than one.

\section{The Crest Factor} \label{SectionCF}

Before briefly summarising more systematically what we mean by the crest factor of 
a solution of a given PDE, let us first introduce some essential notation and functional 
setting \cite{Adams,CoFo,Temam,Robinson}.

\subsection{Basic Definitions and Setup}

Let $\Omega$ be a bounded or unbounded set in $\mathbb{R}^d$, where $d$ is the spatial dimension,
for example the $d$-dimensional torus $\Omega = [0,L]^d$ of side $L$.
For $x\in \Omega$ we write, in general, $x=(x_1,\ldots,x_d)$;
if $d=2$ and $d=3$ in the following, when convenient, we also write $(x_1,x_2)=(x,y)$ and $(x_1,x_2,x_3)=(x,y,z)$, respectively.

For any $p \geq 1$, consider the Banach  
space $L^p(\Omega)$ of continuous functions
$\phi\!: \Omega \times \mathbb{R} \rightarrow \mathbb{R}$, with norm 
$$\|\phi\|_p := \left(\int_{\Omega} |\phi(x,t)|^p \,{\rm d}x \right)^\frac{1}{p}$$ 
and define the sup-norm (or $L^{\infty}$ norm) as
\begin{equation*} 
  \|\phi\|_{\infty} := \sup_{x \in \Omega} |\phi(x,t)| .
\end{equation*}
Here and in the following, the supremum is meant to be taken on the variable $x$ only, and the dependence on time $t$
of the involved norms is not explicitly shown.
For $p=2$, in particular, $L^2(\Omega)$ is the Hilbert space of square integrable functions on $\Omega$,
namely of functions $\phi\!:\Omega \times \mathbb{R} \to\mathbb{R}$ such that $\|\phi\|_{2}<+\infty$.
Given a multi-index $\vec{n}=(n_1,n_2,\ldots,n_d)$, with all the $n_i$ 
non-negative integers, set $|\vec{n}|:=n_1+\ldots+n_d$ and
\begin{equation*} 
D^{\vec{n}} := \frac{\partial^{|\vec{n}|}}
{\partial x_{1}^{n_1} \partial x_{2}^{n_2} \cdot \cdot \cdot \partial x_{d}^{n_d}} ,
\end{equation*}
and consider the Sobolev space of $n$-times differentiable functions on $\Omega,$ 
with up to the $n$-th order derivative in $L^2(\Omega),$ that is 
\begin{equation} \label{Hnze}
H^n(\Omega) := \biggl\{ \phi : \int_{\Omega} (D^{\vec{n}} \phi)^2 {\rm d}x < + \infty
\hbox{ for all } \vec{n} \hbox{ such that } |\vec{n}| = n \biggr\}.
\end{equation}
For vector functions $\phi=(\phi_1,\ldots,\phi_d)\!:\Omega \times \mathbb{R} \to\mathbb{R}^d$, we set
\[
\|\phi\|_p :=\left( \sum_{i=1}^d \int_{\Omega} |\phi_i(x,t)|^p \,{\rm d}x \right)^\frac{1}{p} , \qquad
\|\phi\|_\infty = \max_{i=1,\ldots,d} \|\phi_i \|_{\infty} .
\]

Given a solution $u=(u_1,\ldots,u_d)$ of a PDE in dimension $d$, 
we define the semi-norms 
\begin{equation} \label{Jn}
J_n :=
\!\!\!\!\!\!
\sum_{\substack{n_1,\ldots,n_d \ge 0 \\ n_1+ \ldots + n_d = n}} 
\!\!\!
\frac{n!}{n_1! \ldots  n_d!}  \|D^{\vec{n}} u \|_{2}^2 =
\sum_{i=1}^d \sum_{\substack{n_1,\ldots,n_d \ge 0 \\ n_1+ \ldots + n_d = n}}
\!\!\!
\frac{n!}{n_1! \ldots  n_d!}  \|D^{\vec{n}} u_i \|_{2}^2. 
\end{equation}
Note that, for $n=0$, \eqref{Jn}
reduces to the squared $L^2$ norm $\|u\|_2^2.$
In (\ref{Jn}), provided the solution $u$ is sufficiently regular, which is often the case~\cite{BV,Temam,Robinson},
we naturally identify the functions having the same mixed partial derivatives.

Also, if we consider the $d$-dimensional torus $\Omega = [0,L]^d$ and,
for any infinitely differentiable periodic real function $\phi \! : \Omega\times\mathbb{R} \to \mathbb{R}$,
we take its  Fourier series

\begin{equation} \nonumber 
 \phi (x,t) =\sum_{{k}  \in \mathbb{Z}^d}\, \phi_{{k}} \, e^{2 \pi i {k} \cdot {x} /L},
 \qquad\qquad \overline{\phi}_{{k}}=\phi_{-{k}} ,
\end{equation}
with the Fourier coefficents $\phi_k$ depending on $t$ in general, then, from Parseval's identity, we have
\begin{equation*} 
\|D^{\vec{n}} \phi \|_{2}^2 = L^d \, \left(\frac{2 \pi}{L}\right)^{2n}  \sum_{{k} \in  \mathbb{Z}^d } |{k}|^{2n}|\phi_{{k}}|^2 \,, 
\end{equation*}
with 
$|{k}|^2 = {k} \cdot {k} = k_1^2 + k_2^2 + \ldots + k_d^2$.
By the same token the definition of Sobolev space \eqref{Hnze} can be extended
to any real number $s$ as
\begin{equation*} 
 H^s(\Omega) := \biggl\{  \phi=\sum_{{k} \in \mathbb{Z}^d} \phi_{{k}} \, e^{2 \pi i {k} \cdot {x} /L} :
\overline{\phi}_{{k}}=\phi_{-{k}} \hbox{ and } \sum_{{k} \in  \mathbb{Z}^d} |{k}|^{2s}|\phi_{{k}}|^2 < + \infty \biggr\},
\vspace{-.1cm}
\end{equation*}
and the corresponding norm is given by
\begin{equation*} 
\| \phi \|_{H^s}^2 :=  L^d \, \left(\frac{2 \pi}{L}\right)^{2s}  \sum_{{k} \in  \mathbb{Z}^d } |{k}|^{2s}|\phi_{{k}}|^2\,.
\end{equation*}

Our analysis of the various estimates obtained in this paper involves also taking long-time 
averages of several norms such as the energy or the enstrophy for the Navier-Stokes equations,
where by long-time average we mean the following. Given a continuous function $\phi(t),$ 
which depends upon the solution $u(x,t)$ of a given PDE, with a given initial condition,
for example 
$$ \phi(t) = \int_\Omega (u(x,t))^2 {\rm d}x, $$ 
we define the long-time average of $\phi(t)$ along the solution $u(x,t)$ the quantity \cite{BV,CoFo,Temam}
\[
\langle \phi \rangle := \limsup_{t \rightarrow \infty}\left[\frac{1}{t}  \int_0^t \phi(s) \, {\rm d}s \right].
\]

\subsection{The Crest Factor and Analysis of the Solutions}

Consider a linear or nonlinear PDE in a domain $\Omega \subset \mathbb{R}^d$ and let $u(x,t)$ denote any of its solutions,
with $x \in \Omega$ and $t$ the time. We suppose that the solutions of our PDE possess all the
classical properties such as existence, uniqueness, regularity and continuous dependence on the
initial conditions. The boundary conditions will be specified case by case, including
periodic and Dirichlet boundary conditions. In this setting we define the \emph{crest factor} (CF) 
as the ratio between the sup-norm and the $L^2$  norm of solutions, namely
\begin{equation} \label{CFunb}
C_f \, := \frac{\|u\|_{\infty}}{J_0^{\frac{1}{2}}},
\end{equation}
with $J_0= \|u\|_{2}^2$ and $\|u\|_{\infty} = \max_{i=1,\ldots,d} \| u_i \|_{\infty}.$ 
Thus, the CF of a given solution of a PDE ``captures'' the fluctuation and 
distortion of the amplitude measured by the sup-norm with respect to the $L^2$ norm of that particular solution.
In fact, the CF, as defined in \eqref{CFunb}, is a well-known quantity in electrical engineering and signal analysis,
where it is used, for instance, to measure 
the maximum of the velocity with respect to its average at the points of a grid in wind
tunnels in real turbulence experiments in fluid dynamics \cite{Voke}.

One can see that dimensionally the CF is a length to the power $-d/2$. 
In the case of bounded domains it may be naturally made dimensionless. For example,
by choosing a $d$-dimensional cube of length $L$, namely $\Omega = [0,L]^d$,
we can redefine the CF as
\begin{equation} \label{CFb}
C_f \, := L^\frac{d}{2} \frac{\|u\|_{\infty}}{J_0^{\frac{1}{2}}}.
\end{equation}
and, in such a case, we stick to this definition instead of \eqref{CFunb}.

Let us investigate further our definition of the CF to see what its important features are.
First one can see that, in bounded domains with smooth boundaries, the CF is always bounded below.
In particular, if the domain is the cube of length $L$, then the CF in \eqref{CFb} is larger than  one. Indeed, one has
\[
J_{0}^{\frac{1}{2}} \leq \|u\|_{\infty} L^\frac{d}{2} 
\]
and therefore it follows
that $\displaystyle{1 \leq C_f}$. An upper bound is obtained by using the Gagliardo-Nirenberg
inequality (GNI). Here we take mean zero functions to make the analysis more transparent.
Of course, in general, the mean of the solution is not zero and further terms appear in the formulae below;
for non-zero mean functions see \cite{BaGe}.
Thus, using a GNI for the $L^\infty$ norm of the 
solution we obtain, for any $n > d/2$,
\begin{equation} \label{GNI}
\| u \|_\infty \leq c(n,d,L) \left(\frac{J_n}{J_0}\right)^\frac{d}{4n} J_0^{\frac{1}{2}},
\end{equation}
where $c(n,d,L)$ is a suitable constant depending on $n$, $d$ and $L$,
but not on the solution itself. Eventually, dividing by $J_0^{\frac{1}{2}}$ and 
multiplying by  $L^\frac{d}{2}$, we arrive at bounding \eqref{CFb} by
\[
1 \leq C_f \leq  c(n,d,L) L^{\frac{d}{2}} \left(\frac{J_n}{J_0}\right)^\frac{d}{4n}.
\]

Note that there is a connection between the CF and the length scale $l$ of a given solution,
as defined in \cite{BDGM,DG,GB}, namely 
\begin{equation*} 
l^{- \frac{d}{2}} := c(n,d,L) \left(\frac{J_n}{J_0}\right)^\frac{d}{4n}.
\end{equation*}
Indeed, we find that
\begin{equation} \label{CFscales}
C_f = L^{\frac{d}{2}} \frac{\|u\|_\infty } 
{J_0^{\frac{1}{2}}} \leq L^{\frac{d}{2}} c(n,d,L) \left(\frac{J_n}{J_0}\right)^\frac{d}{4n}
= \left(\frac{L}{l} \right)^{\frac{d}{2}}.  
\end{equation}

Hence effectively by measuring the crest factor of a solution one obtains
some very important features shown by the solution itself. In particular, in bounded domains,
 if the crest factor is of order one then the dynamics is
relatively ``mild'', in the sense that the solution does not have major excursions 
in space-time. By contrast, when the sup-norm of the solution becomes much 
larger with respect to its spatial average, it suggests that the solution does have 
significant fluctuations in space-time; if these fluctuations are chaotic and  intermittent
they represent one of the characteristic signature of hard turbulence. In the light of this 
we will therefore give the following

\begin{criterion*}
A solution of a PDE shows intermittent turbulent behaviour if its 
corresponding CF is a time-dependent chaotic function with strong intermittency.
\end{criterion*}

This definition encapsulates the essential features of intermittent turbulent and 
non-turbulent solutions of PDEs. In the following we will endeavour to elucidate our 
definition of intermittent turbulent behaviour by analysing some classical PDEs and their solutions in 
the various cases where turbulent features are either present or absent.
Various examples have already being  investigated in \cite{MVB,MBPRS,BDG,BaGe}.
In this work we extend that analysis, thereby  strengthening the fact that the CF is a very good tool
for the understanding of some of the  aspects of turbulence.

\section{Linear PDEs} \label{linearPDEs}

In this section we analyse the behaviour  of solutions of some representative linear PDEs
whose solutions can be found by using the method of separation of variables.
We will analyse a non-separable case in Subsection \ref{SSP}.
We first study linear PDEs on a bounded smooth domain with periodic or Dirichlet boundary conditions.

\subsection{The Heat Equation} \label{heat}

We start with the set of linear parabolic separable PDEs. The prototype 
of this class of PDEs is the heat equation on the torus. In space dimension $d=1$,
taking for instance $\Omega=[0,\pi]$ with Dirichlet boundary conditions, we have that  
the function $u=u(x,t)$ of the two variables $(x,t)$ obeys the linear PDE 
\begin{subequations} \label{HE}
\begin{align}
& u_t = k u_{xx}, \hspace{2.8cm} t > 0 ,
\label{HEa} \\
& u(x, 0) = \phi(x),
\label{HEb} \\
& u(\pi,t)= u(0,t)=0 , \hspace{1.0cm} t \ge 0 ,
\label{HEc}
\end{align}
\end{subequations}
where $k$ is a positive constant and $\phi(x)$ is the initial condition, which is a prescribed
integrable function. Here the solution $u(x,t)$  is cast into the form 
\[
u(x,t) = \sum_{n\in\mathbb{Z}} a_n X_n(x) T_n(t),
\]
with the functions
$X_n(x)$ fixed by the boundary conditions, the functions $T_n(t)$ being decaying exponentials and the 
coefficients $a_n$ depending on the initial condition. More precisely, the solution can be  written as
\[
u(x,t) = \sum_{n\in\mathbb{Z}} a_n X_n(x) e^{- k \lambda_n t},
\]
where $X_n$ are the eigenfunctions and $\lambda_n$ are the eigenvalues of the Laplacian operator. 

Here, in addition to the analysis done in \cite{BaGe}, 
we elucidate further a very important feature of solutions of linear PDEs,
which are made of a finite or infinite linear superposition involving all the ``harmonics" of the Fourier series.
Taking any of these ``modes" and analysing its  evolution in space and time, one can see that the corresponding CF is constant.
t is easy to see that the analysis extends to the case where all the modes are involved,
Indeed, assuming, without loss of generality, the eigenfunctions $X_n(x)$ to be normalised so that $\|X_n\|_2=1$ for all $n\in\mathbb{Z}$,
we have, according to \eqref{CFb},
\[
C_f \le \mathfrak{S}_f := \sqrt{\pi} \,
\frac{\displaystyle{\sum_{n\in\mathbb{Z}} a_n \max_{x\in[0,2\pi]} |X_n(x)|
e^{-(\lambda_n-\lambda_0)t}}}{\displaystyle{ \Bigl( \sum_{n\in\mathbb{Z}}  |a_n|^2 e^{-2(\lambda_n-\lambda_0)t} \Bigr)^{\frac12}}} .
\]
For large $t$, we may bound
\[
\mathfrak{S}_f \approx \sqrt{\pi} \max_{x\in[0,2\pi]} |X_0(x)| ,
\]
so that, asymptotically in time, the CF is bounded by a constant which does not depend on the initial condition.
In the case of finite time we refer to Appendix A: the CF is still found to be of order one.

To sum up, according to our definition of turbulent solutions, the solutions of the heat 
equation on the torus cannot have turbulent behaviour.
Here we are using in an essential way  the fact that our PDE is linear, and so
we can use the superposition principle.

In more than one spatial dimension the analysis can be performed in a similar way. For example in two 
dimensions with Dirichlet boundary conditions on a square $\Omega = [0,L]^2$ one has
\begin{subequations} \label{HE2}
\begin{align}
& u_t = k (u_{xx} + u_{yy}), \hspace{1.8cm}(x,y) \in \Omega, \quad\;\;\, t > 0, \\
& u(x,y,0) = \phi(x,y), \hspace{1.66cm} (x,y) \in \Omega, \\
& u(x,y,t)=0 , \hspace{2.67cm} (x,y) \in \partial \Omega, \quad t \ge 0,
\end{align}
\end{subequations}
where the initial condition $\phi(x,y)$ is a prescribed integrable function in his domain. 
In this case the solution can be expressed as
\[
u(x,y,t) = \sum_{n,m=1} a_{n,m} X_n(x) Y_m(y) e^{- k \lambda_{nm} t},
\]
where the functions $\displaystyle{X_n(x)Y_m(y)}$ and the values $\lambda_{nm}$ are the 
eigenfunctions and eigenvalues of the Laplacian in two spatial dimensions for the fixed boundary 
conditions. This class of solutions retain the same structure in any smooth bounded domain $\Omega.$
Hence, by investigating the CF for each Fourier mode, one can see that the analysis above 
in one spatial dimension holds true in any smooth bounded domain as well. So the solutions of 
the heat equation on smooth bounded domains with Dirichlet boundary conditions cannot have 
turbulent behaviour according to our definition, given that the CF for each mode is constant.
The same holds for periodic boundary conditions, with minor changes in the analysis. 

\subsection{The Wave Equation} \label{wave}
 
We now consider the class of linear hyperbolic separable PDEs. The prototype of
this set of equations is the wave equation in $\Omega=[0,L]^d$, for instance in three 
dimensions, i.e. with $d=3$,
\begin{subequations} \label{wave}       
\begin{align}
& u_{tt} = k ( u_{xx} + u_{yy} + u_{zz} ) , \hspace{1.23cm} 
(x,y,z) \in \Omega , \qquad t > 0 , \\
& u(x,y,z,0) = \phi(x,y,z) , \hspace{1.47cm} (x,y,z) \in \Omega , \\
& u_t(x,y,z,0) = \psi(x,y,z) , \hspace{1.47cm} (x,y,z) \in \Omega ,
\end{align}
\end{subequations}
on the torus or on a cuboid with right angles with Dirichlet boundary conditions.
The analysis, \emph{mutatis mutandis}, follows closely that of the heat equation. 
In fact the fundamental solution 
here is made up of an infinite superposition of trigonometric oscillating modes both in 
space and time. So each mode of these solutions will follow a space-time periodic or
quasi-periodic evolution and the corresponding CF is again of order one as for the heat equation.
So, according to our definition of turbulent solutions, the solutions of the wave equation as 
well cannot have turbulent behaviour. 
Similar results are obtained for any other linear hyperbolic separable PDEs on smooth domains.

\subsection{The time-dependent Stokes problem} \label{TDSP}

The time-dependent Stokes problem on a domain $\Omega\subset \mathbb{R}^2$
is obtained as a linearisation of the Navier-Stokes equations, and is described by the equations
\begin{subequations} \label{td-Sokes}       
\begin{align}
& u_{t} = \nu \Delta u - \nabla \mathcal{P} + f ,
\label{N-Sa} \\
& \nabla \cdot u = 0,
\label{N-Sb} \\
& u(x,0)=\phi(x) , \qquad x \in \Omega ,
\label{N-Sd}
\end{align}
\end{subequations}
with $\mathcal{P}$, $\nu$, $f$ and $\phi(x)$ being the pressure, the kinematic viscosity, the external force and the initial condition, respectively.
As shown in \cite{Jiang,Jiangetal}, the solution $u(x,t)$ can be expressed through the Green formula in terms of the fundamental solution,
which in turn is given by the $2\times 2$ matrix $G(x,y,t)$ with entries \cite{GT}
\[
\begin{aligned}
G_{ij} (x,y,t) & = \frac{e^{-\frac{|x-y|^2}{4\nu t}} }{4\pi\nu t} \biggl( \delta_{ij} - \frac{(x_i-y_i)(x_j-y_j)}{|x-y|^2} \biggr) \\
& - \frac{1\!-\!e^{-\frac{|x-y|^2}{4\nu t}} }{4\pi|x-y|^2} \biggl( \delta_{ij} - 2 \frac{(x_i-y_i)(x_j-y_j)}{|x-y|^2} \biggr) ,
\end{aligned}
\]
where $|x-y|=\textstyle{\sqrt{(x_1-y_1)^2+(x_2-y_2)^2}}$ and $\delta_{ij}$ is the Kronecker delta.
Thus, the CF is found to decay in $t$ and no turbulent behaviour appears.

\section{Nonlinear PDEs} \label{nonlinearPDEs}

We now endeavour to analyse the structure of the solutions and the corresponding CF of 
some important nonlinear PDEs. Some cases have already been investigated in \cite{BaGe},
for example the Korteweg-de Vries equation which is the prototype of the so-called 
completely integrable non-linear PDEs having the 
celebrated soliton solution \cite{CalDeg,DEGM}, and the two-dimensional incompressible 
Navier-Stokes equations. Here we extend further our analysis by studying additional 
representative nonlinear PDEs  with the aim to show the power 
that the CF has in extracting the turbulent or non-turbulent nature of their solutions. 

\subsection{Burgers' Equation} \label{SectionBE}

We start with the extensively studied quasi-linear diffusion equation known as Burgers'  
equation \cite{JK}.
We consider the case of a finite interval and, for simplicity,
we choose $\Omega = [0,\pi]$,
with Dirichlet boundary conditions, namely
\begin{subequations} \label{BE}
\begin{align}
& u_{t} = \epsilon \, u_{xx} - u u_x, \qquad\qquad\, 0 \leq x \leq \pi, \\
& u(x,0) = \phi(x),\\
& u(0,t) = u(\pi,t) =0, \hspace{1.1cm} t >0 ,
    \end{align}
\end{subequations} 
with $\epsilon>0$ and $\phi$ being the initial condition.
This equation arises in a number of physical contexts where viscous and nonlinear
effects are equally important and has been introduced as a toy model for mathematical turbulence,
until Hopf and Cole proved that it  can be integrated explicitly \cite{Hopf,Cole}.
Accordingly, here we show that it does not show any turbulent behaviour. In fact, nowadays
the importance of Burgers' equation lies in the modelling of the long-term behaviour of a shock layer \cite{LO};
recently, there has been a renewal of interest in Burgers' equation, including extensions of 
the original one-dimensional model,
with focus on the problem of the so-called Burgers turbulence \cite{BK}.

By using the Cole-Hopf transformation $u \mapsto v$ \cite[\S1.7.1]{JK}, defined as
$$u = -2 \epsilon \frac{v_x}{v},$$
Burgers' equation becomes the heat equation, namely
\begin{subequations} \label{HEbis}
\begin{align}
& v_{t} = \epsilon \, v_{xx}, \hspace{2.9cm} 0 \leq x \leq \pi, \\
& v(x,0) = \alpha \, g(x), \\
& v_x(0,t) = v_x(\pi,t) =0, \hspace{0.74cm} t> 0 ,
\end{align}
\end{subequations}
where $\alpha$ is an arbitrary positive constant and 
\begin{equation} \label{al-g}
g(x) := e^{-\frac{1}{2\epsilon} \int_{0}^{x} \phi(s) \, {\rm d} s }  .
\end{equation}
By using a standard separation of variables method one obtains the solution
\begin{equation} \label{an}
 v(x,t)=\frac{a_0}{2} \!+\! \sum_{n=1}^{\infty} a_n e^{-n^2 \epsilon t} \cos(nx),
\qquad
a_n := \frac{2 \pi}{\alpha} \int_{0}^{\pi} \!\! g(x) \cos(nx)\, {\rm d} x .
\end{equation}
Note that (\ref{al-g}) and (\ref{an}) imply that $a_n > 0$ for all $n\in\mathbb{N}$ and hence $v(x,t) >0$ for all $x\in[0,\pi]$ and all $t\ge0$.
Going back to the function $u(x,t)$ we obtain the solution in the form
\begin{equation} \label{Sol-Bur}
u(x,t)= 2\epsilon \frac{\displaystyle{ \sum_{n=1}^{\infty} n a_n e^{-n^2 \epsilon t} \sin(nx)}}
{\displaystyle{ \frac{a_0}{2} + \sum_{n=1}^{\infty} a_n e^{-n^2 \epsilon t} \cos(nx)} } . 
\end{equation}

We now analyse the CF corresponding to the solution \eqref{Sol-Bur}. The aim is to show that,
according to our definition, the solution does not exhibit any turbulent behaviour. 
First observe that one can view (\ref{Sol-Bur}) as
\begin{equation*} 
u(x,t)= \sum_{n=1}^{\infty} u_n(x,t), \qquad
u_n(x,t) := \frac{\displaystyle{ 2 \epsilon n a_n e^{-n^2 \epsilon t} \sin(nx)}} {\displaystyle{ \frac{a_0}{2} 
+ \sum_{k=1}^{\infty} a_k e^{-k^2 \epsilon t} \cos(kx)} }. 
\end{equation*}

Therefore, for large times, one can write 
\[
u_n(x,t)\approx \frac{4\epsilon n a_n e^{-n^2 \epsilon t} \sin(nx)}{a_0}
\]
and the corresponding CF \eqref{CFb} becomes
\[
C_f \approx \displaystyle{\frac{\sqrt{\pi} \|u_n\|_\infty}{\|u_n\|_2} = \sqrt{2}}.
\]
The case of finite time may be discussed along the lines of Appendix A.

Thus, as in the case of the linear PDEs considered in Section \ref{linearPDEs},
there is not any turbulent behaviour, essentially because the CF for each mode is constant.

From the above analysis it transpires that whenever we are in a context where there 
is separation of the Fourier modes, including asymptotic separation as in the case of 
Burgers' equation, there cannot be any turbulent behaviour.

\subsection{Incompressible Navier-Stokes Equations} \label{SectionINSE}

\noindent We now consider one of the most important systems of classical dynamics, namely
the incompressible Navier-Stokes equations (NSE) on various domains and with various boundary 
conditions.

Let $u$ denote the velocity field, with $u=(u_1,u_2)$ in $d=2$ and \\ 
$u=(u_1,u_2,u_3)$ in $d=3$ with $d$ the spatial dimension; each $u_i = u_i(x,t)$, 
where $x=(x_1,\ldots,x_d)$ and $t$ is the time.
The NSE, with initial condition $u_0(x)$, read
\begin{subequations} \label{N-S}       
\begin{align}
& u_{t} + (u \cdot \nabla) u = \nu \Delta u - \nabla \mathcal{P} + f,
\label{N-Saa} \\
& \nabla \cdot u = 0,
\label{N-Sbb} \\
& u(x,0)=u_0 (x),
\label{N-Sdd}
\end{align}
\end{subequations}
where $\mathcal{P}$ is the pressure, $f$ is the external force applied to the fluid and $\nu$ is 
the kinematic viscosity. Here, as usual, we assume the density of the fluid to be $\rho=1$.

Rigorous results for the Navier-Stokes flow on the two-dimensional flat torus $\Omega=[0 , 2\pi]^2$
show that for any periodic and divergence-free initial condition $u_0 \in J_1$ and
any  time-independent force $f \in L^2(\Omega)$ there is a unique  and global strong solution
which depends continuously on the initial condition $u_0$~\cite{CoFo,BV,Temam,Robinson}.
In $d=3$ it is well known that, with the exception of small initial data or in presence of 
symmetries, we have
global existence of the weak solution and short time existence and uniqueness of the strong 
solution~\cite{CoFo}.
Uniqueness of the weak solution is still an open problem and the global
existence of the strong solution is one of the seven Millennium Prize Problems selected by the Clay Mathematics Institute.

Our analysis of the NSE performed below, in addition to the torus, also investigates
other domains, both bounded and unbounded with corresponding appropriate boundary conditions.
We begin with a famous problem, namely the Stokes second problem.

\subsubsection{Stokes Second Problem} \label{SSP}

Consider an infinitely long plate which is oscillating periodically in its own plane $(y,z)$ 
in the $z$ direction with a periodic motion of the form $z=U_0\cos(\Omega_0 t),$ where $\Omega_0$
and $U_0$ denote the frequency and the amplitude, respectively \cite{Esteban}. The plate has 
on one of his sides an incompressible viscous fluid which extends indefinitely along the 
positive upward $x$-axis, so that it occupies the semi-infinite domain
$V= [0, +\infty) \times ( - \infty, + \infty) \times (- \infty, + \infty) \subset \mathbb{R}^3$.
Furthermore, if gravity is the only external  force acting on the fluid, the only effect of the 
pressure it to compensate its effect \cite{Bat}.
Then, the incompressible Navier-Stokes equation reduces to a linear PDE which describes
the evolution of the third component $u_3$ of the velocity field, which depends only on the 
spatial coordinate $x$;
consistently with the notation introduced above, we are writing 
$(x_1,x_2,x_3)=(x,y,z)$, so that we may set
for simplicity $u_3(x,y,z,t)=w(x,t)$. The linear PDE one obtains is given by \cite{Bat,Esteban}
\begin{subequations} \label{SSPE}
\begin{align}
& w_t = \nu \, w_{xx},
\label{SSPEa} \\
&  w(0,t) = U_0 \Omega_0 \cos(\Omega_0 t),
\label{SSPEb} \\
&  \lim_{x \rightarrow +\infty} w(x,t) = 0 .
\label{SSPEc}
\end{align}
\end{subequations}
This is again the linear diffusion equation which we investigated earlier, that is the heat equation
\eqref{HE} with $k=\nu$. Here, however, besides the special initial condition,
we have different boundary conditions. In fact,
due to the periodic motion of the plate, the separation of variable method does not work.
Therefore we look for a solution in the form 
\begin{equation} \label{solSt}  
  w(x,t) = U e^{\alpha x} e^{i(\beta x + \Omega t)} ,
\end{equation}
for suitable constants $U$, $\Omega$, $\alpha$ and $\beta$ to be fixed,
and then we take its real part. Substituting and then doing the necessary calculations with the
given boundary conditions we find that a solution of the form \eqref{solSt} to \eqref{SSPE} is
\begin{equation} \label{realSt}  
  w(x,t) = U_0 \Omega_0 e^{- \sqrt{\frac{\Omega_0}{2 \nu}} x} 
  \cos \biggl(\Omega_0 t - \sqrt{\frac{\Omega_0}{2 \nu}} x \biggr).
\end{equation}
Looking at (\ref{realSt}) on can see that the solution is in general quasi-periodic in $(x,t),$
when the two frequencies $\Omega_0$ and $\sqrt{\Omega_0/2 \nu}$
are incommensurate, and it becomes periodic when the two frequencies are commensurate. 

Actually, it may be more convenient to investigate the presence of turbulence in terms of 
the vorticity rather than the velocity. The vorticity field ${\bf \omega} := \nabla \times u$ becomes  
\begin{equation} \label{vortic-Stokes}
\begin{aligned}
& {\bf \omega}(x,y,z,t) \!=\! (\partial_{y} u_3 \!-\! \partial_{z} u_2,  \partial_{z} u_1 \!-\! \partial_{x} u_3,
\partial_{x} u_2 \!-\! \partial_{y} u_1) \!=\! (0, -\partial_{x} w, 0) \phantom{\Big)} \\ 
& = \biggl( \! 0, \sqrt{\frac{U_0^2\Omega_0^2}{2\nu}}e^{- \sqrt{\frac{\Omega_0}{2 \nu}} x} 
\biggl( \! \cos\Bigl( \! \Omega_0 t \!-\! \sqrt{\frac{\Omega_0}{2 \nu}} x \! \Bigr) \!-\!
\sin\Bigl( \! \Omega_0 t \!-\! \sqrt{\frac{\Omega_0}{2 \nu}} x \! \Bigr) \! \biggr),0 \! \biggr).
\end{aligned}
\end{equation} 
From (\ref{vortic-Stokes}) we obtain
\begin{equation*} 
\|\omega\|_\infty \!=\! \sqrt{\frac{U_0^2\Omega_0^2}{\nu}}, \qquad
\|\omega\|_2 \!=\! \sqrt{\frac{U_0^2\Omega_0^2}{2\nu}} \left(\frac{2\nu}{\Omega_0}\right)^{\!\!\frac{1}{4}}
\!\! \sqrt{\frac{2+\cos(2\Omega_0t) - \sin(2\Omega_0 t)}{4}} .
\end{equation*}
The corresponding CF \eqref{CFunb} is given by
\begin{equation} \label{CF-vort-S-norms}
C_f = 
\frac{\|\omega\|_\infty}{\|\omega\|_2} 
\sqrt{\frac{4}{2 - \sin(2\omega t)}}
= \left( \frac{2\Omega_0}{\nu}\right)^{\!\!\frac{1}{4}} \!\!
\sqrt{\frac{4}{2+ \cos(2\Omega_0 t)  - \sin(2\Omega_0 t)}} .
\end{equation}

Note that a large amplitude $U_0$ yields large values of the two norms, but does not affect the CF.
We ought to stress that looking at the values of the $\|\omega\|_\infty$ and $\|\omega\|_2$ 
individually, gives only a partial view of the effective turbulent behaviour of the above solution 
of the NSE. Even in the case of large frequency or small viscosity, \eqref{CF-vort-S-norms} shows 
that  the solution exhibits no turbulent behaviour, since the CF is periodic in time.
In this example, the solution decays in space because the exponential damps down quickly 
any complexity in the region not close to the $x$-axis.
In other words the relevant behaviour of the solution is confined in the region close the the
$(y,z)$ plane for $x$ rather small; this region is the so-called \emph{boundary layer}.
However, the analysis above shows that not even inside the boundary layer turbulent behaviour 
may occur. Therefore, the two $L^\infty$ and $L^2$ norms being large individually does not 
characterise turbulent behaviour. 
Besides the example considered above, this is easily seen even just by considering stationary solutions
of the NSE -- or more generally any stationary solution of any partial differential equation.
By a judicious choice of the parameters present in the
given solution one can construct very large $L^\infty$ and $L^2$ norms, but because they are all 
time-independent, there is not, obviously, any turbulent behaviour.
This is clearly illustrated by the example contained in \cite{Kar,Bat}, which we briefly review here.

Consider the stationary solution of the NSE which describes the steady flow of a
viscous fluid in a semi-infinite region bounded by
an infinite disk which rotates in its own plane with angular velocity $\Omega$.
Suppose also that the fluid at infinity has an arbitrary uniform angular velocity $\gamma$ 
about the axis of rotation of the disk.
K\'arm\'an's original problem \cite{Kar} corresponds to $\gamma=0$ and $\Omega \neq 0,$
whereas $\gamma \neq 0$ and $\Omega = 0$  describes a flow which is rotating uniformly
at infinity and which is bounded by a stationary disk.
Under these conditions, even though the flow has not been determined in detail for any
value of the ratio $\gamma/\Omega$ in $(-\infty, +\infty)$,
numerical integrations and analytical computations for some values of the ratio, for instance 
when it is close to 1, show that there exist solutions which are independent of time 
\cite[\S5.5(c)]{Bat} and hence, a fortiori, exhibit no sign of turbulence.

Similar flows which strengthen our point can be found in \cite[\S3.5.1 and \S3.7.1]{RD},
including the work by Stuart which generalises K\'arm\'an's problem \cite{TS} and
the works by Landau and Squire on the round jets \cite{LA,SQ}, which provide further 
examples of time-independent solutions of the NSE.

In the light of the above discussion one can infer that it is the CF which provides the criterion to
extract the real nature of a given solution regarding its  time fluctuations. 
Thus, going back to the CF \eqref{CF-vort-S-norms} for the Stokes second problem,
we can see that we are not in the presence of any real turbulent behaviour.
Note that, if the kinematic viscosity is small or the frequency is large, one can have large 
fluctuations in the evolution of the CF. 
However these fluctuation are periodic and therefore, according to our criterion, we do 
not have genuine turbulent behaviour. In addition, the fluctuations in the CF do not depend 
upon the amplitude of the periodic oscillations. In this kind of situations, the computation 
of the CF automatically rules out any possible turbulent behaviour.

\subsubsection{Motion in a Half-Plane} \label{MHP}

Another instructive example is the motion of a two-dimensional viscous fluid moving on the 
upper half-plane with a zero gradient of pressure \cite{MaPu}.
The NSE \eqref{N-Sa} with $\nabla \mathcal{P} = 0$, in dimension $d=2$ and in the absence of 
the external force, satisfies, besides the continuity equation for incompressible 
fluids \eqref{N-Sb}, the initial conditions
\begin{subequations} \label{MaPu-IC}
\begin{align}
&  u_1(x_1,x_2,0) = U = \hbox{constant} , \qquad x_2 > 0,\\
& u_2(x_1,x_2,0) = 0,
\end{align}
\end{subequations}
and the boundary condition on the ``wall'' $x_2=0$ given by 
\begin{equation} \label{MaPu-BC}
u(x_1,0,t) = 0 .
\end{equation}
The solution is \cite[\S1.4]{MaPu}
\begin{subequations} \label{NS-solution}
\begin{align}
  u_1(x_1,x_2,t) & = \frac{2U}{\sqrt{\pi}} \int_0^\eta e^{-y^2} dy,  \qquad
  \h \eta^2 := \frac{x_2^2}{4\nu t}, \\
  u_2(x_1,x_2,t) & =0. \phantom{\int}
  \end{align}
\end{subequations}
Note that
\[
\begin{aligned}
&  \displaystyle{\lim_{t \rightarrow 0^+} \! u_1(x_1,x_2,t)} \quad \hbox{ for } x_2>0 , \qquad 
\displaystyle{\lim_{x_2 \rightarrow +\infty} \!\!u_1(x_1,x_2,t) =U} , \\
& \displaystyle{\lim_{x_2 \rightarrow 0^+} \! u_1(x_1,x_2,t)} = 
\displaystyle{\lim_{t \rightarrow +\infty} \! u_1(x_1,x_2,t)} = 0.
\end{aligned}
\]

As it is well known \cite{Bat}, the stationary irrotational solution
$(u_1,u_2)=(U,0)$ of the corresponding Euler equation in a half-plane with initial conditions 
\eqref{MaPu-IC} differs significantly from the solution \eqref{NS-solution} of the NSE only 
close to the boundary at $x_2=0$ (within a distance of order $\sqrt{\nu t}$).
This region is the boundary layer where turbulence may normally occur. 

Coming back to the solution \eqref{NS-solution}, the corresponding vorticity field is
\begin{equation} \label{vortic-half}
\omega(x_1,x_2,t) = \partial_{x_1} u_2 (x_1,x_2,t) - \partial_{x_2} u_1 (x_1,x_2,t) =
- \frac{U}{\sqrt{\pi \nu t}} e^{-\frac{ x_2^2 }{4\nu t}}.
\end{equation} 
Note that
$\displaystyle{\lim_{t \rightarrow 0^+} \omega(x_1,x_2,t)} = 
\displaystyle{\lim_{x_2 \rightarrow +\infty} \omega(x_1,x_2,t)} =0$;
it is important to stress that the vorticity is generated by the boundary condition \eqref{MaPu-BC}.
The solution \eqref{vortic-half} is very similar to the solution of the diffusion 
equation on an infinite domain, in the case of a concentrated unit source of heat at the
origin (say) which is switched on for an instant only \cite{JK}. 

We now compute the CF to see
if there are significant fluctuations in the vorticity field, which might suggest that
some turbulence behaviour is present in the dynamics of our problem.
From (\ref{vortic-half}) one can see that 
\begin{equation*} 
   \|\omega\|_\infty = \frac{U}{\sqrt{\pi\nu t}}, \qquad
   \|\omega\|_2 = \frac{U}{(2\pi)^\frac{1}{4}}\frac{1}{(\nu t)^{\frac{1}{4}}}, \qquad t>0 ,
\end{equation*}
and hence 
\begin{equation} \label{CF-Vort}
C_f = \frac{\|\omega\|_\infty}{\|\omega\|_2} = \left(\frac{8}{\pi e^2}\right)^\frac{1}{4}
\frac{1}{(\nu t)^\frac{1}{4}}, \qquad t>0. 
\end{equation}

Therefore, the CF is large only for a limited interval of time and does not depend on the value of $U$.
If, on the one hand, it does not convey the information that the solution is non-trivial only within the boundary layer, namely close to the horizontal $x_1$-axis, on the other hand, because of the decay 
in time, it shows that turbulent behaviour never arises in this case.

\subsubsection{Incompressible Navier-Stokes Equations on the  Torus} \label{SectionINSET}

So far, we have considered situations where an explicit solution can be found. In this section 
we investigate the solutions
of the NSE on the torus in a general setting, where explicit solutions are not known and hence
the time-pointwise analysis of the CF, as well as the analysis of its asymptotic behaviour, 
is extremely difficult -- if not impossible. Therefore, we resort to using the time average 
of our quantities, so as to extract some information about the solution. 

Because of the operation of time averaging, naturally, we lose any information about the time 
evolution. In particular, we cannot detect the 
presence of strong time fluctuations and hence we cannot conclude whether a given 
solution exhibits intermittent turbulence. Nevertheless, the averaged CF  may provide information 
which allows one to exclude the presence of turbulent behaviour,
for instance if it turns out to be of order one.
In this section we will assume that our initial condition and our forcing function
have zero mean, namely
$$\int_\Omega u_0(x) dx =0, \hspace{1cm} \int_\Omega f(x) dx =0.$$
This will allow us to make use of embedding inequalities such as the Gagliardo and Nirenberg 
inequality for mean zero functions. 

\vspace{.4cm}
\noindent\textbf{Two-dimensional case}
\vspace{.2cm}

\noindent In this section we obtain the CF for the two-dimensional incompressible 
NSE on a flat torus $\Omega = [0,2\pi]^2$.
The CF for this case was thoroughly analysed  in \cite{BaGe}. For the sake of the reader
and as a introduction to the three-dimensional case,
we report here the main result. First, we need to define semi-norms similar to \eqref{Jn}, 
which involves the external force, namely 
\begin{equation}\label{forcing}
   \Phi_n := \sum_{i=1}^{d}\sum_{\substack{n_1,\ldots,n_d 
   \ge 0 \\ n_1+ \ldots + n_d =n}}\frac{n!}{n_1! \ldots n_d!} 
      \|D^{\vec{n}} f_i\|_{2}^2,  
\end{equation}
where $f_i$ are the components of the periodic time-independent forcing function.
We wish to add $\Phi_n$ to $J_n$ as defined in \eqref{Jn}, for each $n$, and so in order to make
the involved quantities dimensionally equivalent we multiply each $\Phi_n$ by the
quantity $\tau^2 := L^4 \nu^{-2}$, with $L=2\pi$ in our case
(for more details see \cite{BDGM,DG}). Hence we obtain the quantities
\begin{equation} \label{FN}
F_n := J_n + \tau^2 \Phi_n.
\end{equation}  
Note that $J_n \leq  F_n$ and $\tau^2 \Phi_n \leq F_n.$ 
Then, as discussed in \cite{BDGM,BaGe},
in order to avoid dividing by semi-norms which may become small, we modify once more 
the definition of the CF and set
\begin{equation*} 
C_f \, := L^\frac{d}{2} \frac{\|u\|_{\infty}}{F_0^{\frac{1}{2}}} ,
\end{equation*}
by taking into account also the time-independent, space periodic and mean zero external force.

In \cite{BaGe} we have made some mild restriction
on the structure of the forcing function,
namely we assumed that it has a {\it cut-off} in its spectrum.
This is expressed mathematically by choosing forcing functions such that
they have a smallest length scale \cite{BDGM,DG}
\begin{equation} \label{lambdaf}
\lambda_{f}^{-2} := \sup_{n\in\mathbb{N}} \left[\frac{\Phi_{n+1}}{\Phi_{n}}\right].
\end{equation}

With this setting we have proved the following (for details see \cite{BaGe})

\begin{theorem}
For small values of $\nu$, the long time-averaged crest 
factor for the two-dimensional incompressible NSE on the torus $[0,2\pi]^2$ obeys the estimate
\begin{equation} \label{BGI9}
\langle C_f \rangle
= \biggl\langle L \frac{\|u\|_{\infty}}{F_0^\frac{1}{2}} \biggr\rangle
\leq \,\frac{1}{\sqrt{4\pi}} \frac{L}{\lambda_0} \left[\hat{\eta} + 
\frac{1}{2} \ln \left(\frac{\Phi_0 \Phi_1}{\nu^2} \right) + 
\frac{1}{2e}\frac{1}{\tau^2 \Phi_1} \right]^\frac{1}{2},
\end{equation} 
where
\begin{equation} \label{constants}
\lambda_{0}^{-2} := \lambda_{f}^{-2} + L^{-2}, \qquad  \tau := L^2 \nu^{-1}, 
\qquad 
\hat{\eta} := \eta - \frac{1}{4} \ln(4 c_1 c_2), 
\end{equation}
with $L=2\pi$, $\eta \approx 1.82$, and $c_1$ and $c_2$ 
two positive constants such that $c_1 + c_2 = 1$.
\end{theorem}

Looking at the theorem above one can see that, for very small $\nu,$ the 
crest factor for the solutions of the two-dimensional NSE on the torus, because 
of the term  $\ln  ( \Phi_0 \Phi_1 \nu^{-2})$,
goes like $\nu^{-2}$ logarithmically; this is effectively saying
that, when the viscosity is relatively small, the crest factor of the solutions 
``follows'' the smallest scale of the forcing
function through the term $L/\lambda_0$, whereas, when the viscosity is very small, 
the logarithmic correction becomes relevant.
In other words one  can say that the solutions of the the two-dimensional NSE on 
the torus do not generally show strong turbulent behaviour unless the viscosity parameter 
is extremely small and so very close to the inviscid limit which coincide (on the torus) with the
two-dimensional Euler equations; for the study of this very important problem on various boundary
conditions and various domains see, for example, 
\cite{CoFo,MaPu,CO94,CO2007,CO2008,BeSp,CO2015,CoVi,BaTi19,KVW,KNVW} and the references listed in these papers.
However, in the argument above, small values of viscosity -- and hence large values 
of CF -- does not imply necessarily, in general, hard turbulence since we are not considering 
fluctuations in time (since we are taking the long-time average).

\vspace{.4cm}
\noindent\textbf{Three-dimensional case}
\vspace{.2cm}

\noindent The study of the three-dimensional NSE equations still presents challenging difficulties
because there are not many strong results as in the two-dimensional case. The main open problem is one
does not yet know about existence, uniqueness and regularity of the solution for all time for
large Reynolds numbers: as mentioned earlier in this paper, this open problem is  one of the 
Millennium Prize Problems.
However, one can compute a ``conditional'' CF for the  three-dimensional NSE as a function 
of the sup-norm of the spatial derivatives of the velocity field, namely
$\|Du\|_{\infty}=\max\{ \| \partial u_i/\partial x_j\|_{\infty} : i,j=1,\ldots,d\}$.

First, we rewrite some results from \cite{BDGM,DG} about the time evolution 
for the $F_n$ in (\ref{FN}) in the $d=3$ case and for $n=0,1$. For $n=0$ one has \cite[\S6.5.2]{DG}
\begin{equation} \label{F0}
\frac{1}{2} \dot{F}_0 \le - \nu F_{1} + \nu \lambda_0^{-2} F_0  ,
\end{equation} 
with $\lambda_0$ as in \eqref{constants}.
Dividing by $F_0$ and time averaging we find
\begin{equation} \label{TAF0}
\biggl\langle \frac{F_1}{F_0} \biggr\rangle \leq \lambda_0^{-2}.
\end{equation}
Note that the left hand side term, once integrated, is zero
because $F_0$ is bounded above and below
provided $\tau^2 \Phi_0 > 1.$

The time evolution for $F_1$ is given by \cite[\S6.2]{DG}
\begin{equation} \label{F1}
\frac{1}{2} \dot{F}_1 \le - \nu F_{2} + (c_1 \|D u\|_{\infty}  +\nu \lambda_0^{-2})F_1,
\end{equation}  
where $c_1=c(1,3,L)$, with the notation in \eqref{GNI}.
Actually, the estimate of the constant $c_1$ can be improved with respect to the
value provided by the GNI (see Appendix B).
Likewise to our analysis of $F_0,$ we divide by $F_1$ and time average to obtain 
\begin{equation} \label{TAF1}
\biggl\langle \frac{F_2}{F_1} \biggr\rangle \leq c_1 \, \nu^{-1} \big< \| Du \|_{\infty} \big> + \lambda_0^{-2} ,
\end{equation}
where we have used that the integrated left hand side term in \eqref{F1} is zero as for (\ref{F0})
provided $\tau^2 \Phi_1 > 1.$

Looking at (\ref{TAF1}) we can see that we need solutions having the norm $\|Du\|_{\infty}$ -- 
and hence, \emph{a fortiori}, the norm $\|\omega\|_\infty$ as well - bounded for all time.

Before stating the theorem on the estimate of the CF for the $3d$ NSE on the torus we need the following 

\begin{lemma*}
Given a real number $\alpha$ such that  $0 < \alpha \le 1/2$, 
and given two time-dependent functions $A(t)$ and $B(t)$ having an integrable long-time average,
then the long-time average $\langle (AB)^\alpha \rangle$ satisfies 
\[
\big\langle\bigl(AB\bigr)^\alpha \big\rangle \leq \big\langle A \big\rangle^\alpha  
\big\langle B \big\rangle^\alpha.
\]
\end{lemma*}

\noindent {\it Proof.}
The result is just an application of the H\"older inequality using the conjugate numbers 
$\alpha$ and $1- \alpha.$ Indeed, we have
\[
\big\langle \left(AB\right)^\alpha \big\rangle \leq \big\langle A \big\rangle^\alpha 
\big\langle B^\frac{\alpha}{1-\alpha} \big\rangle^{1-\alpha} \leq 
\big\langle A \big\rangle^\alpha \big\langle B \big\rangle^\alpha ,
\]
where in the second inequality we have used that $0 < \alpha \leq 1/2$ and
the long-time average property $\langle B^\delta \rangle \leq  \langle B \rangle^\delta$, 
which holds for any real number $\delta\in(0,1)$. 

We can now state the following 
\begin{theorem}
The long time-averaged crest 
factor for the three-dimensional incompressible NSE on the torus $[0,L]^3$ obeys the estimate
\begin{equation} \label{CF3d}
\langle C_f \rangle
= \biggl\langle L^\frac{3}{2} \frac{\|u\|_{\infty}}{F_0^\frac{1}{2}}\biggr\rangle \leq
L^\frac{3}{2} c_2
\lambda_0^{- \frac{3}{4}} 
\left(c_1 \, \nu^{-1} \left\langle \|Du \|_{\infty} \right\rangle + \lambda_0^{-2}\right)^\frac{3}{8} ,
\end{equation}
where $c_1$ and $c_2$ are suitable positive constants. One can take $c_1=3$.
\end{theorem}

\noindent {\it Proof.} By computing the CF one obtains
\[
\begin{aligned} \langle C_f \rangle
& = \biggl\langle L^\frac{3}{2} \frac{\|u\|_{\infty}}{F_0^\frac{1}{2}} \biggr\rangle \leq 
L^\frac{3}{2} c_2
\biggl\langle \left(\frac{F_2}{F_0}\right)^\frac{3}{8} \biggr\rangle
\leq
L^\frac{3}{2} c_2
\biggl\langle \left(\frac{F_2}{F_1}\right)^\frac{3}{8}
\left(\frac{F_1}{F_0}\right)^\frac{3}{8} \biggr\rangle \\ 
& \leq
L^\frac{3}{2} c_2
\biggl\langle\frac{F_2}{F_1}\biggr\rangle^\frac{3}{8} \biggl\langle 
\frac{F_1}{F_0}\biggr\rangle^\frac{3}{8}
\leq
L^\frac{3}{2} c_2 \,
\lambda_0^{- \frac{3}{4}} 
\Bigl( c_1 \, \nu^{-1} \big\langle \| Du \|_{\infty} \big\rangle + \lambda_0^{-2} \Bigr)^\frac{3}{8},
\end{aligned}
\]
where we have used the GNI \eqref{GNI} with $n=2$ and $d=3$,
so that $c_2=c(2,d,L)$,
and the outlined property that $J_n \le F_n$ to obtain the first inequality,
the lemma above to obtain the second inequality and, finally, the estimates 
(\ref{TAF0}) and (\ref{TAF1}) to obtain the last inequality.

Note that \eqref{CF3d}, among various things, corroborates the idea that the quantity 
which ultimately controls the solutions
in the case of the three dimensional NSE is the norm $\|Du \|_{\infty}$ -- a result
that appears in the early work by Constantin and Foja\c{s}
(see \cite{CoFo} and references cited therein).
Indeed, \eqref{CF3d} shows that hard turbulence
may appear only in the regime with large fluctuations of $Du$ in space and time,
which imply a large value of $\langle \| Du \|_{\infty}\rangle$,
where one can approximate
\begin{equation} \label{appCF3d}
\langle C_f \rangle
\lesssim L^\frac{3}{2} c_2 \, \lambda_0^{- \frac{3}{4}} 
\left(c_1 \, \nu^{-1} \bigl\langle \| Du \|_{\infty} \bigr\rangle \right)^\frac{3}{8}.
\end{equation}
What is also interesting in \eqref{CF3d} -- and hence in formula \eqref{appCF3d} -- is the power $3/8$:
reaching a state of hard turbulence requires 
``very'' strong and fluctuating vorticity fields because the power $3/8$ is ``not too large".
It is worth stressing that,
in the absence of strong fluctuations of the vorticity, very small values of the viscosity alone
are not sufficient to ensure hard turbulence because the crest factor would be large but nearly constant.

Another interesting point on formula \eqref{CF3d} is the regime when 
$$\displaystyle{\| Du \|_{\infty} \approx L^{-\frac{3}{2}} \| Du \|_2}.$$
In this situation, the crest factor expressed in terms of the vorticity field $\omega$ 
is close to one (see \eqref{SD2} below).
Thus, if turbulence occurs, it may only be homogeneous and isotropic turbulence.
This means that we may only be in the regime of Kolmogorov turbulence:
the crest factor may be chaotic but without strong intermittency. Here one
expects a flow with small spatial length scales similar to Kolmogorov's dissipation 
length \cite[\S3.3]{DG}. This turns out to be the case, as we are about to show.

Suppose that 
\begin{equation} \label{SD}
    \| Du \|_{\infty} \approx L^{-\frac{3}{2}} \| Du \|_2.
\end{equation}
Then, we have
\[
\frac{\|\omega\|_\infty}{\|\omega\|_2}=\frac{\|\omega\|_\infty}{\| Du \|_2} \leq
\frac{\| Du \|_{\infty}}{\| Du \|_2} \approx L^{-\frac{3}{2}} ,
\]
because $\displaystyle{\|\omega\|_2=\| Du \|_2}$ on periodic boundary conditions, 
so that we find approximately 
\begin{equation} \label{SD2}
C_f = L^\frac{3}{2}  \frac{\|\omega\|_\infty}{\|\omega\|_2} \approx L^\frac{3}{2}
\frac{\| Du \|_{\infty}}{\| Du \|_2} \approx 1 .
\end{equation}
Thus, the crest factor is nearly constant  and close to one. 
Let us now obtain Kolmogorov's dissipation length $\lambda_K$
using the approximations above. Recall that \cite[\S3.3]{DG}
\begin{equation} \label{lambdaK}
\lambda_K := C_K \left(\frac{\nu^3}{\varepsilon}\right)^{\frac{1}{4}} ,
\qquad \varepsilon := L^{-3} \nu \, \langle \| Du \|_2^2 \rangle ,
\end{equation}
where $\varepsilon$ is the average energy dissipation rate and $C_K$ is a universal constant.
Following \cite[\S7.4]{DG}, we consider the length scales $l_{n,r}$, for $n,r\in\mathbb{Z}_+$ 
such that $n > r\ge 0$, with
\[
l_{n,r}^{-2} := \biggl\langle \biggl( \frac{F_n}{F_r} \biggr)^{\frac{1}{n-r}} \biggr\rangle .
\]
For $n=2$ and $r=1$, this gives 
\begin{equation} \label{why}
l^{-2} := l_{2,1}^{-2} = \biggl\langle \frac{F_2}{F_1} \biggr\rangle .
\end{equation}
Using the bound \eqref{TAF1}, the definition \eqref{lambdaK} and the approximation \eqref{SD} gives
\begin{equation} \label{SK}
\begin{aligned}
l^{-2} & \le 
c \, \nu^{-1} \big\langle \| Du \|_{\infty} \big\rangle + \lambda_0^{-2} \phantom{\Biggr)} \\
& \le c\,  L^{-\frac{3}{2}} \Biggl(\frac{\nu^{\frac{1}{2}}}{\nu^{\frac{3}{2}}}\Biggr) 
\bigl\langle \| Du \|_2 \bigr\rangle  + \lambda_0^{-2} \leq 
c \, L^{-\frac{3}{2}} \nu^{-\frac{3}{2}} 
\bigl( \nu \bigl\langle \| Du \|_2^2 \bigr\rangle \bigr) ^\frac{1}{2}  + \lambda_0^{-2} \\
& \le c \,\left(\frac{\epsilon}{\nu^3}\right)^\frac{1}{2} + \lambda_0^{-2} =
c \, \left[\left(\frac{\epsilon}{\nu^3}\right)^\frac{1}{4}\right]^2 + \lambda_0^{-2} = 
c _K\, \lambda_K^{-2} + \lambda_0^{-2}, \phantom{\Biggr)}
\end{aligned}
\end{equation}
where we have set $c_K:=c \, C_K^{2}$. Note that, relying on \cite[\S6.2]{DG},
the bounds \eqref{SK} turn out to hold also for the length scales $l_{n,n-1}$, with $n\ge 2$,
which we have not considered explicitly in the discussion above. From \eqref{SK} we conclude that,
when the crest factor is close to 1, Kolmogorov's dissipation length is the natural 
average length scale. Thus, Kolmogorov's dissipation length is the baseline and digressions 
from this are spiky solutions in the dissipation range.

\section{Conclusions} \label{sec5}  

Turbulence -- and intermittent turbulence notably -- is known to be an extremely complicated 
problem. In this work we have stated a criterion for detecting turbulence characterised by 
strong intermittence in
a solution of a particular PDE: \textit{a solution of a PDE shows intermittent turbulent behaviour
if its corresponding CF is a time-dependent chaotic function with strong intermittency.}  

In particular, if the CF of a given solution is time-independent or has a smooth dependence on time,
for instance is time periodic or quasi-periodic, then the solution does not show turbulent behaviour.
According to this, we have checked that the CF of solutions of linear PDEs confirms that they do not
have any turbulent behaviour in their evolution. We have explicitly considered the case of
the heat equation  and of the wave equation,
but the analysis extends to any linear separable PDEs. Another example where the time-dependence
is non-chaotic is provided by the time-dependent Stokes problem.

We have then investigated some nonlinear PDEs. Our first model was Burgers' equation
on a finite interval with Dirichlet boundary conditions. The equation can be explicitly
solved and the CF is easily computed: it shows that, in such a case as well, no turbulent 
behaviour appears (even though Burgers' equation was long considered a prototype for turbulence).

Furthermore, we have considered the 
incompressible Navier-Stokes equations and we have selected some representative
situations which we believe are important to elucidate some features regarding
the presence or absence of intermittent turbulence.
We first considered the Stokes second problem, which can be reduced to a non-separable linear PDE.
The CF turns out to be a periodic function of time and therefore, according to our criterion, 
any kind of turbulent behaviour is absent.
Next, we have analysed the
motion of a two-dimensional viscous fluid moving on the upper half-plane
with a constant horizontal initial velocity field, a zero gradient of pressure
and zero velocity on the ``horizontal wall" for all time. Here we have computed the
CF for the vorticity field and found that it has an algebraic time-decaying 
behaviour, which automatically excludes any turbulent behaviour.

So far we have mentioned situations where we do have an explicit expression of the
solution of the problem under investigation. In all these case we did not find 
any turbulent intermittent behaviour. This scenario therefore suggests that 
turbulence cannot be expressed in closed or explicit form, but it has in its
core a stochastic and erratic nature with many many different length scales 
interacting in a very strong nonlinear way. In the light of this, the last part 
of this paper was devoted to the analysis of the two-dimensional and 
three-dimensional NSE on the torus. Here, apart from relatively easy cases, 
generally it is not possible to provide an explicit solution, and so we have
resorted to taking time averages of various quotients of norms of the solution.

In the two-dimensional case we have found a CF which goes like $\nu^{-2}$ logarithmically. Thus,
for relatively small values of $\nu$ the time averaged CF of the solution 
mimics the smallest scale of the forcing function through the term $L/\lambda_0$.
In the case of very small values of $\nu$ the logarithmic term becomes relevant
giving a large CF. This is effectively saying that in $d=2$ the solutions of the
NSE on the torus can not generally have intermittent turbulent behaviour unless
we are very close to the inviscid limit thereby going to the Euler equations.
In this limit if we are in the presence of strong time fluctuations we can in
principle have intermittent turbulent behaviour.

In the three-dimensional case, for the time averaged CF we found, for large values 
of $\langle \| Du \|_{\infty}\rangle$, the very expressive formula \eqref{appCF3d}. 
We believe this formula is interesting because it implies that, in order
to reach a regime of intermittent turbulence, one must have very strong
fluctuations of spacial variations of the velocity field. 
Since the spatial derivatives of the velocity field generate the vorticity field,
we expect the vorticity field as well to have strong and intermittent turbulent behaviour. 

Of course, when considering the time-averaged CF,
we can no longer detect the presence of time intermittency.
However, computing the averaged CF still allows us to draw some conclusions. Indeed, 
depending on its value, we can conclude what kind of turbulence, if any, is possible:
\begin{enumerate}
\item The CF is of order one ($C_f \approx 1$). This is the regime of soft or mild turbulence,
because we are in a situation  where the $L^\infty$ and the $L^2$ norm of the solution 
are close to each other.
Therefore  no strong intermittent excursions are experienced by the solution.
\item The CF is large ($C_f \gg 1$). This is the regime of hard turbulence.
In this case, if the solution has relevant and intermittent excursion away from its mean 
square space average, then we are in the regime of intermittent turbulence.
\end{enumerate} 

An interesting problem to investigate would be to combine the criterion above with numerical 
analysis of PDEs, including the NSE on the torus. If a numerical solution is available in some 
form and its behaviour is not transparent, then from the computation of the CF one may deduce 
whether the solution is turbulent or not-turbulent.
If no chaoticity in time appears, turbulence must be ruled out. If the CF is chaotic but 
the fluctuations are negligible, then we may have turbulence \textit{\`a la} Kolmogorov. 
Finally, if besides chaoticity, also strong intermittency is present, then we are in the 
regime of intermittent turbulence. \\ 

\noindent {\bf Acknowledgements.}
We thank the anonymous Reviewer for suggesting us to show that, if
one assumes $\displaystyle{\| Du \|_{\infty} \approx L^{-\frac{3}{2}} \| Du \|_2}$,
then one obtains Kolmogorov's dissipation length.
  
\appendix

\section*{Appendix A} \label{A}

With the notations of Section \ref{heat}, set
\[
b_n(t) := |a_n| \,T_n(t), \qquad B(t) := \sum_{n\in\mathbb{Z}} b_n(t) .
\]
For any $\delta>0$, set also $N_1(\delta) := \{ n \in \mathbb{Z} : b_n(t) \ge \delta\}$ and
$N_2(\delta) := \{ n \in \mathbb{Z} : b_n(t) < \delta\}$,
and fix $\delta=\delta(t)$ so that
\[
\sum_{n\in N_1(\delta(t))} b_n(t)\ge \frac{B(t)}{2} .
\]
Then we obtain
\[
B(t) = \sum_{n\in N_1(\delta(t))} b_n(t) 
+ \sum_{n\in N_2(\delta(t))} b_n(t) \le
\sum_{n\in N_1(\delta(t))} b_n(t) + \frac{B(t)}{2}
\]
and hence
\[
B(t) \le 2 \sum_{n\in N_1(\delta(t))} b_n(t)
\le 2 \delta(t) \sum_{n\in N_1(\delta(t))} \frac{b_n(t)}{\delta(t)} \le
 2 \delta(t) \sum_{n\in N_1(\delta(t))} \left( \frac{b_n(t)}{\delta(t)} \right)^2 ,
\]
which implies that 
\[
\frac{\displaystyle{\sum_{n\in\mathbb{Z}} b_n(t)}}{\displaystyle{\sum_n b^2_n(t)}}
\le \frac{1}{\delta(t)} \frac{\displaystyle{\sum_{n\in\mathbb{Z}} 
(b_n(t)/\delta(t))}}{\displaystyle{\sum_{n\in\mathbb{Z}} (b_n(t)/\delta(t))^2}} 
\le \frac{2}{\delta(t)} \frac{\displaystyle{\sum_{n \in N_1(\delta(t))} (b_n(t)/\delta(t))}}{\displaystyle{\sum_{n\in N_1(\delta(t))} (b_n(t)/\delta(t))^2}} 
\le \frac{2}{\delta(t)} .
\]
In conclusion, if we set $K_0= \sup \{ |X_n(x)| : x\in[0,2\pi] , \; n \in \mathbb{Z} \}$, we obtain
\[
\begin{aligned}
C_f & \le \sqrt{2\pi} \, K_0
\frac{\displaystyle{ \sum_{n\in\mathbb{Z}} b_n(t)}}{\displaystyle{\Bigl( \sum_{n\in\mathbb{N}} b^2_n(t) \Bigr)^{\frac12}}} \\
& \le \sqrt{2\pi} \, K_0  \bigl( B(t) \bigr)^{\frac12} 
\frac{ \displaystyle{ \Bigl( \sum_{n\in\mathbb{Z}} b_n(t) \Bigr)^{\frac12}}}{\displaystyle{
\Bigl( \sum_{n\in\mathbb{Z}} b^2_n(t)  \Bigr)^{\frac12} }} 
\le K_f(t) := \sqrt{ \frac{4 K_0^2 \pi B(t)}{\delta(t)}} .
\end{aligned}
\] 
Note that the last bound, if we do not assume any condition on the coefficients $a_n$ 
and the function $T_n(t)$,
is nearly optimal. For instance, if we take $b_n(t)=1/(2N+1)$ for $|n| \le N$ 
and $b_n(t)=0$ for $|n|>N$, then we find
\[
\delta(t) = \frac{1}{2N+1} , 
\qquad B(t) = \sum_{n\in\mathbb{Z}} b_n(t) = 1 ,
\qquad \sum_{n\in\mathbb{Z}} b_n^2(t) = \frac{1}{2N+1} ,
\]
so that $N_1(\delta(t))=\mathbb{Z}$, and this gives $C_f=\sqrt{2\pi/(2N+1)} = 
\sqrt{2K_0^2\pi B(t)/\delta(t)}$.

In general, since the functions $T_n(t)$ decay exponentially, we have that both $B(t)$ 
and $\delta(t)$ tend to $a_0$ when $t$ tends to $+\infty$,
so that the quantity $K_f(t)$ remains bounded. 

\section*{Appendix B} \label{B}

In the case of the NSE with periodic boundary conditions, the bound \eqref{F1} 
can be derived directly from the equations \eqref{N-S}
without using the GNI. Using the definition \eqref{Jn} of $J_1$ and the equation \eqref{N-Sa}, 
and shortening for simplicity
$\partial_i = \partial/\partial x_i$, we obtain
\[
\begin{aligned}
\frac{1}{2} \dot J_1 & = \frac{1}{2} \sum_{i,j=1}^3 \frac{{\rm d}}{{\rm d}t} 
\left\| \partial_j u_i \right\|_2^2 =
\sum_{i,j=1}^3 
\int_\Omega \partial_j u_i \partial_j (u_i)_t \, {\rm d}x\\
& =
\sum_{i,j=1}^3 \int_\Omega \partial_j u_i \partial_j
\biggl( - \biggl( \; \sum_{k=1}^3 u_k \partial_k u_i \biggr) + 
\nu \Delta u_i + f_i - \partial_i p \biggr) {\rm d}x ,
\end{aligned}
\]
with $\Omega=[0,L]^3$. We study separately the four contributions:
\begin{itemize}
\item the first one gives
\[
- \sum_{i,j,k=1}^3  \int_\Omega \partial_j u_i \left( \bigl( \partial_j  u_k \bigr) 
\partial_k u_i  + u_k \partial_j \partial_k u_i  \right) {\rm d}x
\]
where, as far as the second term is concerned, integrating by parts, we find
\[
\int_\Omega \partial_j u_i \left( u_k \partial_j \partial_k u_i  \right) {\rm d}x =
- \!\! \int_\Omega \left( \partial_k \partial_j u_i \right)  u_k \partial_j u_i  \, {\rm d}x
- \!\! \int_\Omega \left( \partial_j u_i \right)^2 \left( \partial_k  u_k \right) {\rm d}x
\]
which implies
\[
\sum_{k=1}^3 \int_\Omega \partial_j u_i \left( u_k \partial_j \partial_k u_i  \right) {\rm d}x =
 - \frac{1}{2} \sum_{k=1}^3 \int_\Omega \left( \partial_j u_i \right)^2 
 \left( \partial_k  u_k \right) {\rm d}x = 0 ,
\]
with the last equality following from \eqref{N-Sb}, whereas the first term is bounded as
\[
\begin{aligned}
& \sum_{i,j,k=1}^3 \left| \int_\Omega \left( \partial_j u_i \right) \left( \partial_j u_k \right) 
\left( \partial_k u_i  \right) {\rm d}x \right| \\
& \qquad \le \max_{j,k} \left| \partial_{j} u_{k} \right| \frac{1}{2}
\sum_{i,j,k=1}^3 \int_\Omega \! \left( \left( \partial_j u_i \right)^2 + 
\left( \partial_k u_i \right)^2 \right) {\rm d}x \\
& \qquad \le \frac{3}{2} \| Du \|_\infty \!\! \int_\Omega \! \biggl( \, \sum_{i,j=1}^3 \!\! 
\left( \partial_j u_i \right)^2 \!+\!\!
\sum_{i,k=1}^3 \!\! \left( \partial_k u_i \right)^2 \! \biggr) {\rm d}x \\
& \qquad = 3 \| Du \|_\infty \!\! \sum_{i,j=1}^3 \! \int_\Omega \!\! 
\left( \partial_j u_i \right)^2 {\rm d}x = 3 \left\| Du \right\|_\infty J_1 ,
\end{aligned}
\]
where we have used the Young inequality;
\item 
the second contribution, by integrating by parts, can be written as
\[
\nu \sum_{i,j,k=1}^3 \int_\Omega \partial_j u_i \partial_j \partial_k^2 u_i
= - \nu \sum_{i,j,k=1}^3 \int_\Omega \left( \partial_k \partial_j u_i \right) 
\left( \partial_j \partial_k2 u_i \right) = - \nu J_2 ,
\]
because of the definition \eqref{Jn} of $J_2$;
\item the third contribution is bounded by using the Cauchy-Schwarz 
inequality, so as to obtain
\[
\left| \sum_{i,j=1}^3 \int_\Omega \partial_j u_i \partial_j f_i \, {\rm d}x \right|
\le \sum_{i,j=1}^3 \biggl( \int_\Omega \left( \partial_j u_i \right)^2 \, {\rm d}x \biggr)^{\frac12}
\biggl( \int_\Omega \left( \partial_j f_i \right)^2 \, {\rm d}x \biggr)^{\frac12} \le \left( \Phi_1 J_1 \right)^{\frac12} ,
\]
where we have used the definition \eqref{forcing} of $\Phi_1$;
\item the fourth contribution vanishes because, again by integrating by parts 
and using \eqref{N-Sb}, one has
\[
\begin{aligned}
\sum_{i,j=1}^3 \int_\Omega \partial_j u_i \partial_j \partial_i p \, {\rm d}x & =
- \sum_{i,j=1}^3 \int_\Omega \left( \partial_j \partial_i u_i \right) \partial_j p \, {\rm d}x \\
& =
- \sum_{j=1}^3 \int_\Omega \partial_j \biggl( \, \sum_{i=1}^3 \partial_i u_i \biggr) \partial_j p \, 
{\rm d}x = 0 .
\end{aligned}
\]
\end{itemize}

By collecting together the results above we arrive at
\[
\frac{1}{2} \dot J_1 \le - \nu J_2 + 3 \left\| Du \right\|_\infty J_1 + 
\left( \Phi_1 J_1 \right)^{\frac12} ,
\]
so that, using \eqref{FN} and the fact that $f$ -- and hence $\Phi_n$ -- does not depend on $t$, 
we obtain
\[
\frac{1}{2} \dot F_1 \le - \nu F_2 + \nu \tau^2 \Phi_2 + 3 
\left\| Du \right\|_\infty F_1 + \tau^{-1} F_1 ,
\]
which, once we take into account the definitions \eqref{lambdaf} 
and \eqref{constants}, gives \eqref{F1} with $c_1=3$.

\end{document}